\documentclass[12pt]{article}

\usepackage{latexsym}
\usepackage{amsfonts}
\usepackage{amsmath}

\newtheorem{theorem}{Theorem}[section]
\newtheorem{lemma}{Lemma}[section]
\newtheorem{proposition}{Proposition}[section]

\begin{document}

\title{An Ornstein-Uhlenbeck process associated to self-normalized sums}

\author{Gopal K. Basak\thanks
{Postal Address:
Stat-Math Unit, Indian Statistical Institute, 203 B. T. Road,
Kolkata 700108, INDIA, Email Address: gkb@isical.ac.in} \\ Indian Statistical Institute
 \and Amites Dasgupta\thanks{Postal Address:
Stat-Math Unit, Indian Statistical Institute, 203 B. T. Road,
Kolkata 700108, INDIA, Email Address: amites@isical.ac.in} \\
Indian Statistical Institute}

\date{}

\maketitle

\begin{abstract}

We consider an Ornstein-Uhleneck (OU) process associated to self-normalised sums in i.i.d. symmetric random variables from the domain of attraction of $N(0, 1)$ distribution. We proved the  self-normalised sums converge to the OU process 
(in $C[0, \infty)$). Importance of this is that the OU process is
a stationary process as opposed to the Brownian motion,
which is a non-stationary distribution (see for example,
the invariance principle proved by Csorgo et al (2003, Ann Probab) for
 self-normalised sums that converges to Brownian motion).
 The proof uses recursive equations similar to those that arise in the area of stochastic approximation and it shows (through examples) that
one can simulate any functionals of any segment of the OU process.
The similar things can be done for any diffusion process as well.

\end{abstract}



\newpage

\section{Introducton}

Self-normalized sums have been studied in the literature, specially for statistical estimation when one doesn't know the variance. For the history of the limit theory and for the problems one faces in the case of stable distributions, an excellent source is Logan et al. \cite{logan}. Later authors have investigated various probabilistic aspects, like the Berry-Esseen bounds, the law of the iterated logarithm, large and moderate deviations etc.  Another important question has been the necessary and sufficient conditions for the limiting distribution to be normal, see Gin\'{e} et al. \cite{gine}.

Recently Cs\"{o}rg\H{o} et al. \cite{csorgo} and 
 Ra\v ckauskas and Suquet \cite{rak_suq}
 have studied Donsker's theorem for self-normalized partial sums processes in zero mean i.i.d. random variables from the domain of attraction of the normal distribution (DAN). 
Thus study of self-normalized sums and related processes 
and their functionals is an important area of research. 
In this article we consider an Ornstein-Uhlenbeck process associated to self-normalized sums in i.i.d. symmetric random variables from DAN. The process arises quite naturally (see (\ref{ke})) as the limit of an approximating sequence of processes involving self-normalized sums. Compared to Cs\"{o}rg\H{o} et al. \cite{csorgo} who obtain a nonstationary process (i.e. Brownian motion), we obtain a stationary process. This however does not make the proofs more involved, in fact apart from two references to Gin\'{e} et al. \cite{gine}, our proofs are quite simple. The techniques may also be useful in studying similar stationary processes associated to self-normalized sums when the i.i.d. symmetric random variables come from other domains of attraction. 

The paper is organised as follows. In section 2 we  give the construction and state the main result, theorem \ref{thm_weakconv}. Section 3 deals with applications. The proofs are given in section 4. The interesting part is the simple  proof of tightness, which uses a recursive equation similar to those that arise in stochastic approximation. In the concluding section 5 we mention a related problem involving  stable distributions.

\section{Notations and the main result}
 
Suppose $\{ X_i \}_1^\infty$ is a sequence of i.i.d. symmetric random variables from the domain of attraction of $N(0, 1)$. The self-normalized sums $Y_j$ are defined as follows: 
\[ S_j = \sum_1^j X_i ,\ V_j = (\sum_1^j X_i^2)^{1/2}, \  Y_j = S_j/V_j ,\]
(under conditions that the denominator is not 0, i.e. $X_i^2 \neq 0$, with probability 1).
The partial sums process related to the above sequence studied by Cs\"{o}rg\H{o} et al. \cite{csorgo} is defined as \[ \Big{\{} \frac{S_{[nt]}}{V_n}, t \geq 0 \Big{\}}, \] whose convergence to the Wiener process is part of their results. 

In contrast we want to employ the self-normalized sum sequence more directly.
Let us fix some notational conventions first. Integer subscripts will denote random variables, integer superscripts will denote processes, processes with real numbers in first brackets will denote the processes at the indicated times. With this convention  
following Basak et al. \cite{basak}, first we define the stochastic processes $Y^n(.)$  as follows:
$Y^n(0)=Y_n$, 
\begin{equation}
 Y^n (\sum_{j=n+1}^{l} b_j^2 ) = Y_{l}, \ \  l \geq n + 1, \label{zi} 
\end{equation}
where $b_{j+1}^2 = E(\frac{X_{j+1}^2}{V_{j+1}^2}) = 1/(j + 1)$. At intermediate points the process is obtained by joining the nearest points linearly 

What the above does will be clear if at this point the reader looks at
equation (\ref{ke}) in section 4, which will also explain what stationary limit one can expect.  For the convenience of the reader and for the examples to be studied in the next section, we repeat some of it here. Since $\max \{ k: \sum_n^k \frac{1}{j+1} \leq t \} \sim [n e^t]$, the limiting behaviour of $Y^n(t)$, following (\ref{zi}), can be seen to be given by the limiting  behaviour of 
\begin{equation} 
Y_{[n e^t]} = \frac{S_{[ne^t]}}{V_{[ne^t]}}, \label{yi} 
\end{equation} 
for $t \geq 0$. However, to state the weak convergence result we have to be more precise. 

Our main result is about the weak limit of the sequence $Y^n$, and is stated as follows: 
\begin{theorem}
$Y^n$ converges weakly in $C[0, 1]$ to the stationary Ornstein-Uhlenbeck process with covariance function $e^{- \frac{1}{2}|t - s|}$.
\label{thm_weakconv}
\end{theorem}
The technique for establishing this is through Prohorov's theorem, i.e. showing convergence of finite dimensional distributions and then proving tightness. This is done in the section 4.

Weak convergence in $C[0, \infty)$ also follows from weak convergence in $C[0, T]$, the proof is similar to that of theorem 2.1 for each $T$. Then one can study other functionals of the sequence $Y^n$ through the same functionals of the Ornstein-Uhlenbeck process. As another application, one can consider the connection
 between
 $\{ Y^n( \cdot)\}$ 
and Brownian motion. This shows that the process $\{ \sqrt{t} Y^n(\log t): t \geq 1 \}$ converges weakly to Brownian motion $\{ B(t): t \geq 1 \}$, and one can use Brownian functionals to study functionals of this process and vice-versa.
 Applications of these are considered in the next section.

\section{Applications}

We now present a few examples connecting functionals of Ornstein-Uhlenbeck processes
with those of
self-normalized random walks where
 we can apply our result. This would also indicate the
differences between our result and the invariance principle
presented in Cs\"org\H o et al \cite{csorgo}.

\bigskip
{\bf Example 1.}
For any fixed, $0 \leq s \leq t$,

$$
\mbox{  (a)   } \ \max_{[m e^s] \leq i \leq [m e^t]} Y_{i}
\to \max_{s \leq u \leq t} Y(u) \ \mbox{ in distribution as  } m \to \infty ,
$$
and 
$$
\mbox{  (b)   } \ \min_{[m e^s] \leq i \leq [m e^t]} Y_{i}
\to \min_{s \leq u \leq t} Y(u) \ \mbox{  in distribution as  } m \to \infty ,
$$
where
$Y_{i}=\frac{S_i}{V_{i}}$ and
 $Y(\cdot)$ is a stationary Ornstein-Uhlenbeck process with parameter $(-1/2,1)$,
i.e., it satisfies the stochastic differential equation,
$$
dY(u) = \alpha Y(u) du + \sigma dB(u) , \  \mbox {  with  } \alpha = -(1/2), \ \sigma = 1,
$$
and $\{B(\cdot)\}$ is a standard Brownian motion independent of $Y(0)$, and $Y(0)$ has
$N(0, 1)$ distribution.

\bigskip

{\bf Example 2.}
For any fixed, $0 \le s \le t$,

$$
\mbox{  (a)   } \ \frac{1}{[me^t] - [me^s]}\sum_{[m e^s] < i \le [m e^t]} Y_{i}
\to \int_{s}^{t} Y(u) du \ \mbox{ in distribution as  } m \to \infty ,
$$
$$
\mbox{  (b)   } \ \frac{1}{[me^t] - [me^s]}\sum_{[m e^s] < i \le [m e^t]} |Y_{i}| \to \int_{s}^{t} |Y(u)| du \ \mbox{  in distribution as  } m \to \infty ,
$$
and
$$
\mbox{  (c)   } \ \frac{1}{[me^t] - [me^s]}\sum_{[m e^s] < i \le [m e^t]} Y_{i}^2 \to \int_{s}^{t} (Y(u))^2 du \ \mbox{  in distribution as  } m \to \infty ,
$$
where $\{Y(u)\}$ is a stationary Ornstein-Uhlenbeck process with parameter $(-1/2,1)$, as above.

\bigskip

The above examples can be seen in contrast with the
 invariance principle
presented in Cs\"org\H o et al \cite{csorgo} as follows:

\bigskip

{\bf Example} $\mathbf{1^\prime}$.
For any fixed 
$0 \le t$,

$$
\mbox{  (a)   } \ \max_{1 \le i \le [m e^t]} \frac{S_i}{V_{[me^t]}} \to
 \max_{0 \le u \le 1}
 B(u) \ \mbox{ in distribution as  } m \to \infty ,
$$
and 
$$
\mbox{  (b)   } \ \min_{1 \le i \le [m e^t]} \frac{S_i}{V_{[me^t]}} \to 
\min_{0 \le u \le 1}
 B(u) \ \mbox{  in distribution as  } m \to \infty ,
$$
where
 $\{B(\cdot)\}$ is a standard Brownian motion.

\bigskip

{\bf Example} $\mathbf{2^\prime}$.
For any fixed
 $0 \le t$,

$$
\mbox{  (a)   } \ \frac{1}{[me^t]}\sum_{1 < i \le [m e^t]} \frac{S_i}{V_{[me^t]}} \to 
 \int_{0}^{1} B(u) du
 \ \mbox{ in distribution as  } m \to \infty ,
$$
$$
\mbox{  (b)   } \ \frac{1}{[me^t]}\sum_{1 < i \le [m e^t]} \frac{|S_i|}{V_{[me^t]}} \to 
\int_{0}^{1} |B(u)| du 
\ \mbox{  in distribution as  } m \to \infty .
$$
and
$$
\mbox{  (c)   } \ \frac{1}{[m e^t]}\sum_{1 \le i \le [m e^t]} \Big(\frac{S_i}{V_{[me^t]}}\Big)^2 \to
\int_{0}^{1} (B(u))^2 du 
\ \mbox{  in distribution as  } m \to \infty ,
$$
where
 $\{B(\cdot)\}$ is a standard Brownian motion.

\section{Proofs}
 
{\bf Proof of theorem \ref{thm_weakconv}}

To show that $Y^n$ converges to $Y$ we show convergence of finite dimensional distributions and tightness of $Y^n$. We state these as two propositions. 

In the proof of {\em convergence of finite dimensional distributions} and their identification we use Lemma 3.2 of Gin\'{e} et al. \cite{gine} which adapted to our context reads as: for some slowly varying function $l(n)$ we have 
\begin{eqnarray}
 \frac{S_n}{\sqrt{n l(n)}} &\Rightarrow& N(0, 1), \nonumber \\
  \frac{V_n^2}{n l(n)} &\rightarrow& 1 \mbox{ in probability (Gin\'{e}  et  al.  \cite{gine}) }. \label{re} 
\end{eqnarray}

Gin\'{e} et al. \cite{gine} have stated their lemma for the infinite variance case. In the case of finite variance we can take $l(n) = E X^2, \forall n$, and CLT and SLLN gives the same result. In the following proposition we make no distinction between the finite and infinite variance case.  

\begin{proposition}
The finite dimensional distributions of $Y^n$ converge to those of a stationary Ornstein-Uhlenbeck process with covariance function $e^{- \frac{1}{2}|t - s|}$.
\label{prop_findim}
\end{proposition}
{\bf Proof:} Since $\max \{ k: \sum_n^k \frac{1}{j+1} \leq t \} \sim [n e^t]$, we are interested in the limiting distribution of 
\begin{equation}
(Y_{[ne^{t_1}]}, \cdots , Y_{[ne^{t_k}]})\label{ke}
\end{equation}
 for $0 < t_1 < \cdots
< t_k$ (the above in the case of finite variance is clearly similar to 
\[ \Big(\frac{S_{[ne^{t_1}]}}{\sqrt{[n e^{t_1}]}}, \cdots , \frac{S_{[ne^{t_k}]}}{\sqrt{[n e^{t_k}]}}\Big),\] and this motivates the limiting distribution). It is clear that by independence, and by the slowly varying nature of $l(n)$ 
\[  \Big{(} \frac{S_{[ne^{t_1}]}}{V_{[ne^{t_1}]}}, \frac{S_{[ne^{t_2}]} - S_{[ne^{t_1}]} }{\sqrt{V_{[ne^{t_2}]}^2 - V_{[ne^{t_1}]}^2}}, \cdots , \frac{S_{[ne^{t_k}]} - S_{[ne^{t_{k-1}}]} }{\sqrt{V_{[ne^{t_k}]}^2 - V_{[ne^{t_{k -1}}]}^2}} \Big{)} \] converges in distribution to a vector of i.i.d. $N(0, 1)$, and also that \[ \Big{(} \frac{V_{[ne^{t_1}]}}{V_{[ne^{t_i}]}}, \frac{ \sqrt{V_{[ne^{t_2}]}^2 - V_{[ne^{t_{1}}]}^2}}{V_{[n e^{t_i}]}}, \cdots ,
\frac{ \sqrt{V_{[ne^{t_i}]}^2 - V_{[ne^{t_{i -1}}]}^2}}{V_{[n e^{t_i}]}} \Big{)} \] converges in probability to \[ \Big(\frac{e^{t_1/2}}{e^{t_i/2}}, \frac{\sqrt{e^{t_2} - e^{t_1}}}{e^{t_i/2}}, \cdots , 
\frac{\sqrt{e^{t_i} - e^{t_{i -1}}}}{e^{t_i/2}}\Big), \] for $i \leq k$ (this can be seen by writing \[ V^2_{[ne^{t_i}]} = V^2_{[ne^{t_1}]} + \sum_2^i (V_{[ne^{t_l}]}^2 - V_{[ne^{t_{l -1}}]}^2) ,\] and then using independence after dividing by a $V_{[ne^{t_r}]}^2 - V_{[ne^{t_{r -1}}]}^2$ where $1 < r \leq i$). Hence, writing 
\[ \frac{S_{[ne^{t_i}]}}{V_{[ne^{t_i}]}} = \frac{V_{[ne^{t_1}]}}{V_{[ne^{t_i}]}}\frac{S_{[ne^{t_1}]}}{V_{[ne^{t_1}]}} + \sum_2^i \frac{ \sqrt{V_{[ne^{t_l}]}^2 - V_{[ne^{t_{l-1}}]}^2}}{V_{[n e^{t_i}]}} \frac{S_{[ne^{t_l}]} - S_{[ne^{t_{l-1}}]} }{\sqrt{V_{[ne^{t_l}]}^2 - V_{[ne^{t_{l-1}}]}^2}}\]
and using Slutsky's theorem, the limiting distribution 
of \[ \Big( \frac{S_{[ne^{t_1}]}}{V_{[ne^{t_1}]}}, \cdots , \frac{S_{[ne^{t_k}]}}{V_{[ne^{t_k}]}}\Big)\] is seen to be multivariate normal with  covariance between the limits of  $\frac{S_{[ne^{t_i}]}}{V_{[ne^{t_i}]}}$ and $\frac{S_{[ne^{t_j}]}}{V_{[ne^{t_j}]}}$ 
(taking $i < j$) being \[ \frac{e^{t_1}}{e^{t_i/2} e^{t_j/2}} + \sum_{m = 2}^i \frac{e^{t_m} - e^{t_{m - 1}}}{e^{t_i/2} e^{t_j/2}} = e^{\frac{1}{2}(t_i - t_j)},\] 
completing the proof. \hfill $\Box$

\vspace{.2in}

In the proof of {\em tightness} we need the following two auxiliary results. The first one is part of  theorem 3.3 from Gin\'{e} et al. \cite{gine} stated in a form suitable for our purpose. 
\begin{lemma} (Gin\'{e} et al. \cite{gine}) For $X$ in DAN with $E X = 0$, we have \[ E \Big(\frac{X_{j+1}^4}{V_{j+1}^4}\Big) = o(\frac{1}{j+1}).\]
\hfill $\Box$
\label{lem_gine}
\end{lemma}

\vspace{.2in}

The other result is a simple calculation and we state it as a lemma:
\begin{lemma}
For $j < k$ we have \[ E \Big(\frac{X_{j+1}^2}{V_{j+1}^2} \frac{X_{k+1}^2}{V_{k+1}^2}\Big)
= \frac{1}{(j+1)(k+1)}.\]
\label{lem_crosspdct}
\end{lemma}

{\bf Proof:} We write $V^2_{k+1} = V^2_{j+1} + \widetilde{V}^2_{k - j}$ where 
$\widetilde{V}_{k-j}^2 = \sum_{i = j+2}^{k+1} X_i^2$. We then have, using the i.i.d. nature of $X$'s and the symmetry of the expressions,
\begin{eqnarray*}
1 &=& E \Big(\frac{V_{j+1}^2 (V^2_{j+1} + \widetilde{V}^2_{k - j})}{V^2_{j+1} (V^2_{j+1} + \widetilde{V}^2_{k - j})}\Big) \\
&=& E \Big(\frac{V_{j+1}^2}{V_{k+1}^2}\Big) + \sum_{l = 1, m = j+2}^{j+1, k+1} E \Big(\frac{X^2_l X_m^2}{V_{j+1}^2 (V_{j+1}^2 + \widetilde{V}^2_{k - j})}\Big)  \\
&=& \frac{j+1}{k+1} + (k-j)(j+1)  E \Big(\frac{X_{j+1}^2}{V_{j+1}^2} \frac{X_{k+1}^2}{V_{k+1}^2}\Big) ,
\end{eqnarray*}
from which the result follows. \hfill $\Box$

\vspace{.2in}

Now we prove

\begin{proposition}
The sequence of stochastic processes $Y^n$ is tight.
\label{prop_tight}
\end{proposition}
{\bf Proof:} From results in the literature (or (\ref{re})) $Y^n(0) = Y_n$ is tight (converges to $N(0, 1)$). Thus, by the corollary to  theorem 7.4 of
 Billingsley (\cite{bill}) we only need to show that $\forall \epsilon, \forall \eta, \exists n_0$ and $0 < \delta < 1$ such that \[ \frac{1}{\delta} P \Big{\{} \sup_{t \leq s \leq t + \delta} |Y^n(s) - Y^n(t)| \geq \epsilon \Big{\}} \leq \eta , \ n \geq n_0, \ \forall t \in [0, 1].\] In our problem  it suffices to show that 
 $\forall \epsilon, \forall \eta, \exists n_0$ and $0 < \delta < 1$ such that  
\[ \frac{1}{n + 1} + \cdots + \frac{1}{l + 1} < \delta \Rightarrow P \Big{\{} \sup_{ n \leq k \leq l + 1 } | Y_{k}  - Y_{n} | \geq \epsilon \Big{\}} \leq \eta \delta, \ \forall n \geq n_0. \] 
Note that $Y_j$'s satisfy the following recursion:
\[ Y_{j + 1} - Y_j = - Y_j \frac{X_{j+1}^2}{V_{j + 1}(V_j + V_{j + 1})} + \frac{X_{j + 1}}{V_{j + 1}}.\] Hence \[ |Y_k - Y_{n}| \leq \sum_{n}^{k - 1} \frac{X_{j+1}^2}{V_{j+1}(V_j + V_{j + 1})} |Y_j| + \Big{|} \sum_{n}^{k - 1} \frac{X_{j + 1}}{V_{j + 1}} \Big{|},\]
and 
\begin{equation}
 \sup_{ n  \leq k \leq l + 1 } | Y_{k}  - Y_{n} | \leq \sum_{n}^{l} \frac{X_{j+1}^2}{V_{j+1}(V_j + V_{j + 1})} |Y_j| + \sup_{ n  \leq k \leq l } \Big{|} \sum_{n}^k \frac{X_{j + 1}}{V_{j + 1}} \Big{|}.\label{ai}
 \end{equation} 

We thus have 
\begin{eqnarray}
P \Big{\{} \sup_{ n \leq k \leq l + 1 } | Y_{k}  - Y_{n} | \geq \epsilon \Big{\}} &\leq& 
P\Big(\sum_{n}^{l} \frac{X_{j+1}^2}{V_{j+1}(V_j + V_{j + 1})} |Y_j| > \epsilon/2\Big) \nonumber \\
&+& P\Big(\sup_{ n  \leq k \leq l } \Big{|} \sum_{n}^k \frac{X_{j + 1}}{V_{j + 1}} \Big{|} > \epsilon/2\Big) . \label{ai1}
\end{eqnarray}
Inside the absolute values of the second term of (\ref{ai1}) we have a martingale sequence because the distributions of the $X$'s are symmetric about 0 and we can apply Doob's maximal inequality for submartingales (and then Burkholder's inequality) to get for an appropriate constant $C_4$,  
 \begin{eqnarray}
 && P \Big{\{} 
\sup_{ n \leq k \leq l } \Big{|} \sum_{n}^k
\frac{X_{j + 1} }{V_{j + 1}} \Big{|} > \epsilon /2 \Big{\}} \nonumber \\
 &\leq& C_4 \frac{16}{\epsilon^4}
E \Big( \sum_{n}^{l} \frac{X_{j + 1}^2}{V_{j + 1}^2}\Big)^2 \nonumber \\
&=& C_4
 \frac{16}{\epsilon^4} \Big{\{}
\sum_{n }^{l } a_{j+1} \frac{1}{j+1} + 
\sum_{\begin{array}{ll}
 & j, k = n \\
 & j \neq k 
\end{array}}^l 
\frac{1}{(j+1)(k+1)} \Big{\}}, \label{bi}
\end{eqnarray}
where we have used lemmas \ref{lem_gine} and \ref{lem_crosspdct} and we have written \[ E \Big(\frac{X_{j+1}^4}{V_{j+1}^4}\Big) = a_{j+1} \frac{1}{j+1}.\] Recalling from lemma \ref{lem_gine} that  $a_{j+1} \rightarrow 0$ we have from (\ref{bi})   
\begin{eqnarray}
&& P \Big{\{} 
\sup_{ n \leq k \leq l } \Big{|} \sum_{n}^k
\frac{X_{j + 1}}{V_{j + 1}} \Big{|} > \epsilon /2 \Big{\}} \nonumber \\
&\leq& C_4 \frac{16}{\epsilon^4} \Big{\{} \max_{n \leq j \leq l} \Big|a_{j+1} - \frac{1}{j+1}\Big| 
\Big(\sum_{n }^{l} \frac{1}{j+1}\Big) + \Big(\sum_{n }^{l} \frac{1}{j+1}\Big)^2 \Big{\}} \nonumber \\ 
&=& C_4 \frac{16}{\epsilon^4} \Big{\{} \max_{n \leq j \leq l} \Big|a_{j+1} - \frac{1}{j+1}\Big| 
+ \Big(\sum_{n }^{l} \frac{1}{j+1}\Big) \Big{\}} \Big(\sum_{n }^{l} \frac{1}{j+1}\Big) \label{ci}
\end{eqnarray} 
the first part of which can be made smaller than $\eta$ by choosing $n_0$ large and $\delta$ small remembering  $\sum_{n }^{l} \frac{1}{j+1} < \delta$.

To the first term of (\ref{ai1}) we shall apply Markov's inequality with second moment.
The $X_i$'s being symmetric the distributions remain same if we replace them by $\epsilon_i X_i$ where $\epsilon_i$'s are i.i.d. Rademacher random variables independent of $X_i$'s.
Thus distributionally $|Y_j|$ equals \begin{equation} \frac{| \epsilon_1 X_1 + \cdots + \epsilon_j X_j|}{V_j}.\label{dhi} \end{equation}
 Now the first term of (\ref{ai1}) gives 
\begin{eqnarray}
&& P \Big{\{} \sum_{n}^{l } \frac{X_{j+1}^2}{V_{j+1}(V_j + V_{j + 1})} |Y_j| > \epsilon /2 \Big{\}} \nonumber \\
 &\leq& \frac{4}{\epsilon^2} E \Big( \sum_{n}^{l } \frac{X_{j+1}^2}{V_{j+1}^2} |Y_j| \Big)^2 \nonumber \\
&=& \frac{4}{\epsilon^2} \Big{\{} \sum_{n}^{l } E \Big(\frac{X_{j+1}^4}{V_{j+1}^4} |Y_j|^2\Big) \nonumber \\
&+& 
\sum_{\begin{array}{ll}
 & j, k = n \\
 & j \neq k 
\end{array}}^l 
 E \Big(\frac{X_{j+1}^2}{V_{j+1}^2} \frac{X_{k+1}^2}{V_{k+1}^2} |Y_j|  |Y_k| \Big) \Big{\}}.\label{di}   
\end{eqnarray}

Remembering the expression (\ref{dhi}) taking conditional expectation given \newline $\{ X_1, \cdots , X_{j+1} \}$ we get 
\begin{equation}
 E \Big{\{} \frac{X_{j+1}^4}{V_{j+1}^4} E \Big( |Y_j|^2 \ \Big| X_1, \cdots , X_{j+1} \Big) \Big{\}} 
= E \Big(\frac{X_{j+1}^4}{V_{j+1}^4}\Big) , \label{ei}  
\end{equation}
and supposing $j < k$ applying conditional Cauchy inequality  given \newline $\{ X_1, \cdots , X_{k+1} \}$ we get 
\begin{equation}
E \Big{\{} \frac{X_{j+1}^2}{V_{j+1}^2} \frac{X_{k+1}^2}{V_{k+1}^2} E \Big( |Y_j| \ |Y_k| \ \Big|  X_1, \cdots , X_{k+1} \Big) \Big{\}} 
\leq E \Big(\frac{X_{j+1}^2}{V_{j+1}^2} \frac{X_{k+1}^2}{V_{k+1}^2}\Big) . \label{fi}
\end{equation}
Now
 combining (\ref{ei}) and (\ref{fi}) for (\ref{di}) we have as before 
\begin{eqnarray}
&& P \Big{\{} \sum_{n}^{l } \frac{X_{j+1}^2}{V_{j+1}(V_j + V_{j + 1})} |Y_j| > \epsilon /2 \Big{\}} \nonumber \\
&\leq& \frac{4}{\epsilon^2} \Big{\{} \max_{n \leq j \leq l}\Big|a_{j+1} - \frac{1}{j+1}\Big| 
+ \Big(\sum_{n }^{l} \frac{1}{j+1}\Big) \Big{\}} \Big(\sum_{n }^{l} \frac{1}{j+1}\Big). \label{gi}
\end{eqnarray}
(\ref{ci}) and (\ref{gi}) combined for (\ref{ai1}) give us the required tightness. \hfill $\Box$

\vspace{.2in}

Using propositions \ref{prop_findim} and \ref{prop_tight} the proof of 
theorem \ref{thm_weakconv} is complete.

\section{Concluding remarks}

In the case of symmetric $\alpha$-stable distributions (or symmetric distributions from the domain of attraction of symmetric $\alpha$-stable distributions) it is known from Logan et al. \cite{logan} that \[ \frac{ X_1^4 + \cdots + X_n^4}{V_n^4} \] converges weakly to a positive random variable. From this it follows that \[ n E \Big(\frac{X_n^4}{V_n^4}\Big)\] converges to a positive constant, 
say $c_{\alpha}$. In fact, since (Logan et al. \cite{logan}) $E (S_n^4/V_n^4) \rightarrow 1 + \alpha$ and obviously $E (V_n^2)^2/V_n^4 = 1$, following the proof of theorem 3.3 of Gin\'{e} et al. \cite{gine} and using the simpler symmetry assumption of ours, the constant $c_{\alpha}$ can be shown to be equal to $1 - (\alpha/2)$. Thus our proof of proposition 
\ref{prop_tight} cannot be repeated along the same lines in this case although the necessary calculations for the finite dimensional distributions (corresponding characteristic functions) can be repeated following Logan et al. \cite{logan}.

\end{document}